\documentclass{amsart}
\usepackage{amsmath}
\usepackage{amscd}
\usepackage{amssymb}

\newcommand{\mrm}{\mathrm}

\renewcommand{\div}{\operatorname{div}}

\newcommand{\Vol}{\operatorname{Vol}}

\newcommand{\tr}{\operatorname{tr}}

\newtheorem{theorem}{Theorem}[section]
\newtheorem{theorem/definition}{Theorem/Definition}[section]
\newtheorem{proposition}{Proposition}[section]
\newtheorem{lemma}{Lemma}[section]

\theoremstyle{remark}
\newtheorem{remark}{Remark}[section]

\theoremstyle{definition}

\begin{document}
\title
{On second variation of Perelman's Ricci shrinker entropy}
\author{Huai-Dong Cao and Meng Zhu}
\address{Department of Mathematics\\ Lehigh University\\
Bethlehem, PA 18015, USA} \email{huc2@lehigh.edu \&
mez206@lehigh.edu}

\date{}

\begin{abstract} 
In this paper we provide a detailed proof of the second variation
formula, essentially due to Richard Hamilton, Tom  Ilmanen and the
first author, for Perelman's $\nu$-entropy.  In particular, we
correct an error in the stability operator stated in Theorem 6.3
of  \cite{Cao08b}. Moreover, we obtain a necessary condition for
linearly stable shrinkers in terms of the least eigenvalue and its
multiplicity of certain Lichnerowicz  type operator associated to
the second variation.
 \end{abstract}
\maketitle

\footnotetext[1]{The first author was partially supported by NSF
Grants  DMS-0506084 and DMS-0909581; the second author was
partially supported by NSF Grant DMS-0354621 and a Graduate Summer
Fellowship awarded by the School of Arts and Sciences at Lehigh
University.}

\section{The Results} \label{sec:1}
\hspace{.75cm}

A complete Riemannian metric $g_{ij}$ on a smooth manifold $M^n$
is called a {\it gradient shrinking Ricci soliton} if there exists
a smooth function $f$ on $M^n$ such that the Ricci tensor $R_{ij}$
of the metric $g_{ij}$ satisfies the equation
$$R_{ij}+\nabla_i\nabla_jf=\frac 1 {2\tau} g_{ij} \eqno(1.1)$$
for some constant $\tau>0$. The function $f$ is called a {\it
potential function} of the Ricci soliton.  When $f$ is a constant
we obtain an Einstein metric of positive scalar curvature. Thus,
Ricci solitons include Einstein metrics as a special case.

Ricci solitons correspond to self-similar solutions of Hamilton's
Ricci flow, and often arise as limits of dilations of
singularities in the Ricci flow. In particular shrinking solitons
are possible Type I singularity models in the Ricci flow. We refer
the readers to \cite{Cao08b}, \cite{Cao09} and the references therein for
more information on Ricci solitons.

Ricci solitons can be viewed as fixed points of the Ricci flow, as 
a dynamical system, on the space of Riemannian metrics modulo
diffeomorphisms and scalings.
 In \cite{P1}, Perelman introduced the $\mathcal{W}$-functional
$$
\mathcal{W}(g_{ij},f,\tau)=\int_M [\tau (R+|\nabla f|^2)+f-n](4\pi
\tau)^{-\frac{n}{2}}e^{-f}dV,
$$
on a compact manifold $M^n$, where $g_{ij}$ is a Riemannian
metric, R the scalar curvature, $f$ a smooth function on $M^n$,
and $\tau$ a positive scale parameter. The associated
$\nu$-entropy is defined by
$$
\nu(g_{ij})=\inf\{\mathcal W(g,f,\tau): f\in C^\infty(M), \tau>0,
(4\pi\tau)^{-\frac{n}{2}}\int e^{-f}dV=1\}.
$$
It turns out that the $\nu$-entropy is monotone increasing under
the Ricci flow, and its critical points are precisely given by
gradient shrinking solitons. In particular, it follows that all
compact shrinking Ricci solitons are gradient shrinking solitons,
a fact shown by Perelman \cite{P1}.

In dimensions 2 and 3, Hamilton \cite{Ha88} and Ivey \cite{Iv}
respectively showed that the only compact shrinking solitons are
quotients of the round spheres. However, for dimension $n\ge  4$,
compact non-Einstein shrinking solitons do exist. Specifically in
dimension $n=4$,  Koiso \cite{Ko} and the first author
\cite{Cao94} independently constructed a gradient K\"ahler-Ricci
shrinking soliton on $\Bbb CP^2\#(-\Bbb CP^2)$, and Wang-Zhu
\cite{WZ} on   $\Bbb CP^2\#(-2\Bbb CP^2)$, while in the noncompact
case Feldman-Ilmanen-Knopf \cite{FIK} constructed the
$U(2)$-invariant gradient shrinking K\"ahler-Ricci solitons on the
tautological line bundle $\mathcal{O}(-1)$ of $\Bbb CP^{1}$, the
blow-up of $\Bbb C^{2}$ at the origin. These are the only known
examples of nontrivial (i.e., non-Einstein or non-product)
complete shrinking Ricci solitons  in dimension $4$ so far.

In \cite{CHI04}, Hamilton, Ilmanen and the first author initiated
the study of linear stability of Ricci solitons. They found the second variation
formula of the  $\nu$-entropy for positive Einstein manifolds and investigated
the linear stability of certain Einstein manifolds.
By definition, a Ricci shrinker or Einstein manifold is called {\sl linearly stable}
if the second variation is non-positive. They showed that, while the round sphere
$\Bbb{S}^ n$ and the complex projective space $\Bbb{C}P^n$  are linearly stable,
many known Einstein manifolds are unstable for the Ricci flow so that generic
perturbations acquire higher $\nu$-entropy and thus can never
return near the original metric. In particular,  all K\"ahler-Einstein
manifolds with Hodge number $h^{1,1}>1$ are unstable.

In dimension $n\ge 4$, so far no one knows how to classify
Einstein manifolds of positive scalar curvature, let alone
gradient shrinking Ricci solitons. However, as far as applications
of the Ricci flow to topology is concerned, one is more interested
in stable shrinking solitons since unstable ones could be
perturbed away thus may not represent generic singularities. For
this reason, it is desirable and important to classify stable
shrinking Ricci solitons. Note that, the work of
Cao-Hamilton-Ilmanen \cite{CHI04} suggests that most gradient
shrinking Ricci solitons are unstable. In fact, Hamilton conjectured
that, at least in dimension $n=4$,  compact linearly stable
shrinkers are rank one symmetric spaces, namely either the round
sphere $\Bbb{S}^ 4$ or the complex projective space $\Bbb{C}P^2$
with the Fubini-Study metric.  Of course, in studying linear
stability of shrinkers, the second variation formula of the
$\nu$-entropy is indispensable. In this paper, we present a
detailed proof of the second variation formula, first due to
Hamilton, Ilmanen and the first author  (cf. Theorem 6.3 in
\cite{Cao08b}), for Ricci shrinkers.

To state the second variation formula, we need some notations
first. For any symmetric 2-tensor $h=h_{ij}$ and 1-form
$\omega=\omega_i$, we denote

$$Rm(h,\cdot):=R_{ijkl}h_{jl},$$
$$\text{div}\:\omega:=\nabla_i\omega_i, \qquad (\text{div}\:
h)_i:=\nabla_jh_{ji}.$$ Moreover, as done in \cite{Cao08b}, we
define
$$\div_f \omega:=e^{f}\div( e^{-f}\omega)= \nabla_i\omega_i-\omega_i\nabla_if\eqno(1.2)$$ 
and
$$\div_f h:=e^{f}\div( e^{-f}h)=\div h-h (\nabla f, \cdot), \eqno(1.3)$$ 
i.e.,
$$(\div_f\: h)_i=\nabla_jh_{ij}-h_{ij}\nabla_jf.$$
We also define $\div_f^{\dagger}$ on 1-forms (and similarly on functions) by \\
$$(\div_f^{\dagger} \omega)_{ij}=-(\nabla_i\omega_j+\nabla_j\omega_i)/2
=-(1/2)L_{\omega^\#}g_{ij}\eqno (1.4)$$ so that 
$$\int_M e^{-f} <\div_f^{\dagger}\omega, h>dV= \int_M e^{-f}<\omega, \div_f h>dV.\eqno(1.5)$$
Here $\omega^\#$ is the vector field dual to $\omega$. Clearly,
$\div_f^{\dagger}$ is just the adjoint of $\div_f$ with respect to
the weighted $L^2$-inner product
$$(\cdot, \cdot)_f=\int_M <\cdot, \cdot> e^{-f}dV. \eqno(1.6)$$

Finally we denote 
$$\Delta_f := \Delta -\nabla f \cdot \nabla \eqno(1.7)$$

\begin{remark} If we denote by $\div^{*}$ the adjoint of $\div$ with respect to the $L^2$-inner product
$$(\cdot, \cdot)=\int_M <\cdot, \cdot>dV,$$ then, as pointed out in \cite{Cao08b}, one can easily verify that
$$\div_f^{\dagger}=\div^{*}. \eqno(1.8)$$  
\end{remark}

 Now we can state the full second variation formula for Ricci shrinkers:

\begin{theorem} \label{thm1.1}{\bf (Cao-Hamilton-Ilmanen)} Let $(M^n, g_{ij}, f)$ be a
compact Ricci shrinker with the potential function $f$ and
satisfying the Ricci soliton equation (1.1). For any symmetric
2-tensor $h=h_{ij}$, consider variations
$g_{ij}(s)=g_{ij}+sh_{ij}$. Then the second variation
$\delta^2_g\nu(h,h)$ is given by
\begin{align*}
\left.\frac{d^2}{ds^2}\right|_{s=0}\nu(g(s))
&=\frac{\tau}{(4\pi\tau)^{n/2}}\int_M  <\hat{N}h, h> e^{-f} dV,\\
\end{align*}
where the stability operator $\hat N$ is given by
$$\hat Nh:=\frac{1}{2}\Delta_f h+Rm(h,\cdot)+\div_f^{\dagger}\div_f h+
\frac{1}{2}\nabla^2\hat
v_h-Rc \ \frac{\int_M <Rc, h>e^{-f}}{\int_M Re^{-f}}, \eqno(1.9)
$$
and $\hat v_h$ is the unique solution of
\begin{align*}
\Delta_f \hat v_h+\frac{\hat v_h}{2\tau}=\div_f\div_f
h,\qquad\int_M \hat v_h e^{-f}=0.
\end{align*}
\end{theorem}

\medskip
\begin{remark}\label{rmk1.1} As we pointed out before, Theorem \ref{thm1.1} is essentially due to Hamilton,
Ilmanen and the first author  (cf. Theorem 6.3 in \cite{Cao08b}).
However, the coefficient of the  last term of the stability
operator  $\hat N$ (which depends on $\delta\tau$, the first
variation of the parameter $\tau$)  was stated incorrectly in
\cite{Cao08b}.  One of our contributions in this paper is to
derive an explicit formula for $\delta\tau$ (see Lemma 2.4 below),
thus obtaining the correct coefficient and hence a complete second
variation formula for Ricci shrinkers.  Of course, it would be interesting to investigate the noncompact case as well. In this case, the asymptotic estimates on potential functions and volume growth upper bound proved by Cao-Zhou \cite{CZhou08}, and an integral bound on the Ricci curvature by Munteanu-Sesum \cite{MS} should be very helpful.  
We point out that, while
the stability operator $\hat N$ is already quite useful even
without knowing the explicit coefficient of the last term, it will
be rather crucial to have this explicit and correct coefficient in
efforts of trying to classify stable shrinkers. For example, this
explicit coefficient is essential in showing that the Ricci tensor
is a null eigen-tensor of $\hat N$ (see Lemma \ref{lem3.3}) which rules out
any hope of using the Ricci tensor as a possible unstable
direction.
\end{remark}

\begin{remark}\label{rmk1.2}
In the very recent work \cite{HM}, Stuart Hall and Thomas Murphy
proved that K\"ahler-Ricci shrinking solitons with Hodge number
$h^{1,1}>1$ are unstable, thus extending  the results of
Cao-Hamilton-Ilmanen \cite{CHI04} in the K\"ahler-Einstein case
mentioned above. In the course of their proof, they also verified
the second variation formula stated in \cite{Cao08b}, though
didn't find out explicitly the coefficient of the last term of
$\hat N$ (which does not affect the proof of their result since
they only considered certain special variations orthogonal to
$Rc$).
\end{remark}

\begin{remark} \label{rmk1.3}If $(M^n, g_{ij})$ is Einstein with $Rc=\frac{1}{2\tau} g_{ij}$, Theorem \ref{thm1.1} reduces to

\begin{theorem} \label{thm1.2}{\bf (Cao-Hamilton-Ilmanen \cite{CHI04})}
Let $(M^n, g_{ij})$ be a Einstein manifold and consider variations $g_{ij}(s)=g_{ij}+sh_{ij}$.
Then the second variation $\delta^2_g\nu(h,h)$ is given by
$$\left.\frac{d^2}{ds^2}\right|_{s=0}\nu(g(s))=\frac{\tau}{\Vol (M, g)}\int_M <Nh,h>dV,$$
where
$$
Nh:=\frac{1}{2}\Delta h+Rm(h,\cdot) +\div^*\div\: h+
\frac{1}{2}{\nabla}^2v_h -\frac{g}{2n\tau\Vol(M, g)}\int_M \tr_g h
\: dV,
$$
and $v_h$ is the unique solution of
\begin{align*}
\Delta v_h+\frac{v_h}{2\tau}=\div\div\: h, \qquad\int_M v_h=0.
\end{align*}
\end{theorem}
\end{remark}

\medskip
Finally, using the second variation formula, we obtain the
following necessary condition for linearly stable shrinkers:

\begin{theorem}\label{thm1.3} Suppose $(M^n, g_{ij}, f)$ is a compact linearly stable  shrinking soliton satisfying (1.1),
then $-\frac{1}{2\tau}$ is the only negative eigenvalue of the
operator $\mathcal{L}_f$ (with Rc being an eigen-tensor), defined by
$$\mathcal{L}_f h=\frac{1}{2}\Delta h+Rm(h,\cdot),  \eqno(1.10) $$
on $\ker\div_f$ and the multiplicity of  $-\frac{1}{2\tau}$  is one.  In particular, $-\frac{1}{2\tau}$
is the least eigenvalue of $\mathcal{L}_f $ on $\ker\div_f$.
\end{theorem}

\begin{remark} In proving Theorem \ref{thm1.3}, the explicit coefficient of Rc term in $\hat N$ is not needed.
\end{remark}

\begin{remark} In the mean curvature flow, Colding and Minicozzi \cite{CM} have shown that for any shrinker its mean curvature $H$ is an eigenfunction
of certain operator involved in the corresponding stability
operator, and that for any (linearly) stable shrinker the mean
curvature function $H$ belongs to the least eigenvalue of the
operator which in turn implies that $H$ does not change sign. This
fact and a prior theorem of Huisken allow them to classify compact
stable mean curvature shrinkers. Our Theorem \ref{thm1.3} above can be
considered as the Ricci flow analogy of their results.
\end{remark}

\noindent {\bf Acknowledgements.} The first author would like to
thank Qiang Chen, Richard Hamilton, Tom Ilmanen for stimulating
discussions, and Stuart Hall for helpful communications.

\section{The Proof of Theorem \ref{thm1.1}} \label{sec:2}

In this section, we describe the first variation of the
$\nu$-entropy and derive the second variation formula as stated in
Theorem \ref{thm1.1}.

On any given compact manifold $M^n$, Perelman \cite{P1} introduced
the $\mathcal{W}$-functional
$$
\mathcal{W}(g_{ij},f,\tau)=\int_M [\tau (R+|\nabla f|^2)+f-n](4\pi
\tau)^{-\frac{n}{2}}e^{-f},
$$
where $g_{ij}$ is a Riemannian metric, R the scalar curvature, $f$
a smooth function on $M^n$, and $\tau$ a positive scale parameter.
Clearly the functional $\mathcal{W}$ is invariant under
simultaneous scaling of $\tau$ and $g_{ij}$, and invariant under
diffeomorphisms. Namely, for any positive number $a$ and any
diffeomorphism $\varphi$ we have
$$
\mathcal{W}(a\varphi^*g_{ij},\varphi^*f,a
\tau)=\mathcal{W}(g_{ij},f,\tau).
$$

\begin{lemma} \label{lem2.1}
{\bf (Perelman \cite{P1}, see also Lemma 1.5.7 in \cite{CZ05})}
If $h_{ij}=\delta g_{ij},\; \phi=\delta f,\;\mbox{and}\;
\eta=\delta\tau$, then
$$\arraycolsep=1.5pt\begin{array}{rcl}& &\delta \mathcal{W}(h_{ij},\phi,\eta)\\[4mm]&=& (4\pi\tau)^{-\frac{n}{2}}(\int_M-\tau
h_{ij}(R_{ij}+\nabla_i\nabla_jf-\frac{1}{2\tau}g_{ij})e^{-f}\\[4mm]
&&+\int_M(\frac{1}{2}\tr_g
h-\phi-\frac{n}{2\tau}\eta)[\tau(R+2\Delta f
-|\nabla f|^2)+f-n-1]e^{-f}\\[4mm]
&&+\int_M \eta(R+|\nabla f|^2-\frac{n}{2\tau})e^{-f}).
\end{array}
$$
\end{lemma} 

\medskip
Now, recall that the associated $\nu$-energy is defined by
$$
\nu(g_{ij})=\inf\{\mathcal W(g,f,\tau): f\in C^\infty(M), \tau>0
\},
$$
subject to the constraint
$$(4\pi\tau)^{-\frac{n}{2}}\int e^{-f}=1.\eqno(2.1)$$
One checks that $\nu(g_{ij})$ is realized by
a pair $(f,\tau)$ that solve the equations
$$\tau(-2\Delta f+|Df|^2-R)-f+n+\nu=0,\eqno(2.2)$$
and 
$$(4\pi\tau)^{-\frac{n}{2}}\int fe^{-f}=\frac{n}{2}+\nu.
\eqno(2.3)$$

For any symmetric 2-tensor $h=h_{ij}$, consider variations
$g_{ij}(s)=g_{ij}+sh_{ij}$. Using Lemma \ref{lem2.1},  (2.2) and (2.3), one
obtains the following first variation for the $\nu$-entropy.

\begin{lemma} \label{lem2.2} The first variation $\delta_g\nu(h)$ of the $\nu$-entropy  is given by
\begin{align*}
\frac{d}{ds}\nu(g_{ij}(s)) = & (4\pi\tau)^{-\frac{n}{2}}\int
-\tau <h, Rc+\nabla^2 f-\frac{1}{2\tau} g> e^{-f} dV\\
= & (4\pi\tau)^{-\frac{n}{2}}\int -\tau
h_{ij}(R_{ij}+\nabla_i\nabla_j f-\frac{1}{2\tau} g_{ij})e^{-f} dV.
\end{align*}
\end{lemma}

A stationary point of $\nu$ thus satisfies the Ricci soliton
equation (1.1):
$$R_{ij}+\nabla_i\nabla_j f-\frac{1}{2\tau}g_{ij}=0,
$$
which says that $g_{ij}$ is a gradient shrinking Ricci soliton.

Note that, by diffeomorphism invariance of $\nu$, $\delta_g\nu(h)$
vanishes on Lie derivatives, hence on $h_{ij}=\nabla_i\nabla_j
f=\frac{1}{2}L_{\nabla f}g_{ij}$. By scale invariance it also
vanishes on multiplies of the metric. Inserting
$h_{ij}=-2(R_{ij}+\nabla_i\nabla_j f-\frac{1}{2\tau}g_{ij})$, one
recovers Perelman's  formula that finds that $\nu(g_{ij}(t))$ is
monotone increasing on the Ricci flow, and constant if and only if
$g_{ij}(t)$ is a gradient shrinking Ricci soliton. In particular,
it follows that any compact shrinking Ricci soliton is necessarily
a gradient soliton,
a result first shown by Perelman \cite{P1}.  \\

Now we are going  to derive the second variation formula. \\

\noindent {\it Proof of Theorem 1.1.} \ From the first variation formula in Lemma \ref{lem2.2}, we see that the second variation at a gradient shrinker $(M^n, g_{ij}, f)$ is given by
\begin{align*}
 \delta^2 \nu_g (h, h)= & (4\pi\tau)^{-\frac{n}{2}}\int
-\tau <h, \delta (Rc+\nabla^2 f-\frac{1}{2\tau} g)> e^{-f}\\
= & (4\pi\tau)^{-\frac{n}{2}}\int
-\tau <h, \delta Rc+\delta\nabla^2 f-\frac{1}{2\tau} h> e^{-f} \\
& + (4\pi\tau)^{-\frac{n}{2}}(-\frac{\delta\tau}{2\tau}) \int_M
\tr_g h  e^{-f}.
\end{align*}

\begin{lemma}\label{lem2.3} We have
\begin{align*}
\delta Rc+\delta\nabla^2 f-\frac{1}{2\tau} h = -\frac{1}{2}
\Delta_f h-Rm (h, \cdot) -\div_f^{\dagger}\div_f h-\nabla^2
(-\delta f +\frac{1}{2} \tr_g h) .
\end{align*}
\end{lemma}

\begin{proof} First of all, it is well-known that the variation $\delta Rc$ of the Ricci tensor is given by
$$(\delta Rc)_{ij}= -R_{ikjl}h_{kl} + \frac{1}{2}(\nabla_i\nabla_k h_{jk}  + \nabla_j\nabla_k h_{ik} + R_{ik}h_{jk} + R_{jk}h_{ik} - \Delta h_{ij} - \nabla_i\nabla_j \tr_gh),\eqno(2.4)$$
and, by direct computations (see, e.g., \cite{MZhu}),
$$ (\delta \nabla^2 f)_{ij}=\nabla_i\nabla_j (\delta f)- \frac{1}{2}(\nabla_i h_{jk} + \nabla_j h_{ik} - \nabla_k h_{ij})\nabla_k f.\eqno(2.5)$$ 
On the other hand, by the definition of $\div_f$ and
$\div_f^{\dagger}$ and using the shrinking soliton equation (1.1),
we have
\begin{align*}
\div_f^{\dagger}\div_f h=& -\frac{1}{2}[\nabla_i(\div_f h)_j+\nabla_j(\div_f h)_i]\\
=&  -\frac{1}{2}[\nabla_i(\nabla_k h_{jk}-h_{jk} \nabla_k f) + \nabla_j(\nabla_kh_{ik}-h_{ik} \nabla_k f)]\\
=&  -\frac{1}{2} (\nabla_i\nabla_k h_{jk} + \nabla_j\nabla_k h_{ik}- \nabla_k f\nabla_i h_{jk}  -\nabla_k f\nabla_j h_{ik})\\
& -\frac 12(R_{ik}h_{kj}+ R_{jk}h_{ki}) +\frac {1}{2\tau} h_{ij}.
\end{align*}
Now, combining the above computations, we arrive at
\begin{align*}
\delta Rc+\delta\nabla^2 f = & -\frac{1}{2} \Delta_f h-Rm (h, \cdot) -\div_f^{\dagger}\div_f h\\
& -\nabla^2 (-\delta f +\frac{1}{2} \tr_g h) +\frac{1}{2\tau} h.
\end{align*} 

\end{proof}

Next we derive the variation $\delta\tau$ of the parameter $\tau$.

\begin{lemma} \label{lem2.4}
We have
\begin{align*}
\delta\tau = \tau \frac{\int_M <Rc, h> e^{-f}}{\int_M R e^{-f}}.
\end{align*}
\end{lemma}

\begin{proof} First of all, from (1.1) we get 
$$ R+\Delta f =\frac{n}{2\tau}. \eqno(2.6)$$
Also,  it is  well-known that
$$R+|\nabla f|^2=\frac {f-\nu} {\tau}.\eqno (2.7) $$ 
From (2.6) and (2.7) it follows that
$$-\Delta_f f=:|\nabla f|^2-\Delta f= \frac {f-\nu-n/2} {\tau}.\eqno(2.8)$$
Moreover, from (2.4) and (2.5) and using (1.1), we get 
$$\delta R=-\frac{1}{2\tau} \tr_g h + h_{ij}\nabla_i\nabla_j f + \nabla_i\nabla_j h_{ij} - \Delta \tr_g h, \eqno(2.9)$$
and
$$\delta(\Delta f)= \Delta(\delta f)  - h_{ij}\nabla_i\nabla_j f - \nabla_i h_{ij}\nabla_j f + \frac{1}{2}\nabla_i \tr_g h \nabla_i f \eqno(2.10)$$ respectively.  Also,
$$\delta |\nabla f|^2= 2 \nabla _i f\nabla_j (\delta f)-h_{ij}\nabla_i f\nabla_j f. \eqno(2.11)$$

When we integrate (2.2) against the measure $(4\pi\tau)^{-\frac{n}{2}}e^{-f}dV$ and use (2.3), we obtain
$$(4\pi\tau)^{-\frac{n}{2}}\int_M \tau(|\nabla f|^2 + R)e^{-f} \mrm{d}V = \frac{n}{2}. \eqno(2.12)$$
On the other hand, by differentiating (2.1) and (2.3),  we have
$$(4\pi\tau)^{-\frac{n}{2}}\int_M (-\frac{n}{2\tau} \delta\tau - \delta f + \frac{1}{2} \tr_g h)e^{-f}  =0, \eqno(2.13)$$
and
$$(4\pi\tau)^{-\frac{n}{2}}\int_M f(-\frac{n}{2\tau} \delta\tau - \delta f + \frac{1}{2} \tr_g h)e^{-f}  + (4\pi\tau)^{-\frac{n}{2}}\int_M \delta f e^{-f} =0.\eqno(2.14)$$
Now, differentiating (2.2) and using (2.6), (2.10) and (2.11), we obtain
\begin{align*} 
0 = & \delta\tau (-\frac{n}{2\tau} + |\nabla f|^2 - \Delta f) - \delta f\\
& + \tau (-2\Delta (\delta f)+ 2h_{ij}\nabla_i\nabla_j f + 2\nabla_i h_{ij}\nabla_j f - \nabla_i (\tr_g h) \nabla_i f \\
& + 2\nabla_i f \nabla_i (\delta f) - h_{ij}\nabla_i f \nabla_j f - \delta R).
\end{align*}
Substituting  (1.1) and (2.9) in the above identity, we get
\begin{align*}
0 = &  -\frac{n}{2\tau}\delta\tau - 2\tau \Delta (\delta f) + 2\tau \nabla(\delta f)\nabla f - \delta f + \delta\tau (|\nabla f|^2-\Delta f)\\
&+ \tau(2h_{ij}\nabla_i\nabla_j f + 2\nabla_i h_{ij}\nabla_j f - \nabla_i (\tr_g h) \nabla_i f - h_{ij}\nabla_i f\nabla_j f)\\
& + \tau(\frac{1}{2\tau}\tr_g h - h_{ij}\nabla_i\nabla_j f -
\nabla_i\nabla_j h_{ij} + \Delta \tr_g h).
\end{align*}
But, by definition of $\div_f$, we compute that
\begin{align*}
\div_f \div_f h= & \nabla_i(\nabla_j h_{ij}-h_{ij}\nabla_j f)-\nabla_i f (\nabla_j h_{ij}-h_{ij}\nabla_j f)\\
=& \nabla_i\nabla_j h_{ij}-h_{ij}\nabla_i\nabla_j f-2 \nabla_if
\nabla_j h_{ij}+h_{ij}\nabla_i f\nabla_j f.
\end{align*}
Hence, we get
\begin{equation*}
0=(-\frac{n \delta\tau}{2\tau} - \delta f + \frac{1}{2} \tr_g h) + \delta\tau(-\Delta_f f) +\tau \Delta_f (-2 \delta f + \tr_g h) - \tau \div_f \div_f h.
\end{equation*}

Multiplying the above identity by $f$ and integrating against the measure
$(4\pi\tau)^{-\frac{n}{2}}e^{-f}dV$, we get
\begin{align*}
0 &= (4\pi\tau)^{-\frac{n}{2}}\int_M f(-\frac{n}{2\tau} \delta\tau- \delta f + \frac{1}{2} \tr_g h) e^{-f}dV \\
& \ \ + (4\pi\tau)^{-\frac{n}{2}} \delta\tau \int_M f (-\Delta_f f)e^{-f}\mrm{d}V\\
&\ \  + (4\pi\tau)^{-\frac{n}{2}} \int_M \tau f \Delta_f (-2 \delta f + \tr_g h)e^{-f}dV\\
&\ \  - (4\pi\tau)^{-\frac{n}{2}} \int_M \tau f (\div_f \div_f h) e^{-f}dV.\\
\end{align*}
By (2.14) and integration by parts, the above identity becomes
\begin{align*}
0 &= (4\pi\tau)^{-\frac{n}{2}}\int_M -\delta f e^{-f}\mrm{d}V + \delta\tau (4\pi\tau)^{-\frac{n}{2}}\int_M |\nabla f|^2 e^{-f}\mrm{d}V\\
& \ \ + (4\pi\tau)^{-\frac{n}{2}}\int_M \tau(-2\delta f+\tr_g
h)\Delta_f f e^{-f}\mrm{d}V\\
&\ \ - (4\pi\tau)^{-\frac{n}{2}}\int_M\tau
<h, \nabla^2 f> e^{-f}\mrm{d}V.
\end{align*}
Using (1.1), (2.8) and (2.12), we obtain
\begin{align*}
0 &= -(4\pi\tau)^{-\frac{n}{2}}\int_M \delta f e^{-f}\mrm{d}V + \frac{n}{2\tau}\delta\tau-\delta\tau (4\pi\tau)^{-\frac{n}{2}}\int_M R e^{-f}\mrm{d}V\\
& \ \ + (4\pi\tau)^{-\frac{n}{2}}\int_M 2\tau(\frac{n}{2\tau}\delta\tau +\delta f - \frac{1}{2} \tr_g h)(\frac{1}{\tau}f - \frac{\nu}{\tau} - \frac{n}{2\tau}) e^{-f}\mrm{d}V\\
& \ \ + (4\pi\tau)^{-\frac{n}{2}}\int_M (-\frac{1}{2}\tr_g h +
\tau h_{ij}R_{ij})  e^{-f}\mrm{d}V.
\end{align*}
By using (2.13) and (2.14), we arrive at

\begin{align*}
0 &= (4\pi\tau)^{-\frac{n}{2}}\int_M (\frac{n}{2\tau}\delta\tau + \delta f - \frac{1}{2} \tr_g h) e^{-f}\mrm{d}V - \delta\tau (4\pi\tau)^{-\frac{n}{2}}\int_M R e^{-f}\mrm{d}V\\
& \ \  + (4\pi\tau)^{-\frac{n}{2}}\int_M \tau R_{ij}h_{ij} e^{-f}\mrm{d}V\\
&=  - \delta\tau (4\pi\tau)^{-\frac{n}{2}}\int_M R e^{-f}\mrm{d}V
+ (4\pi\tau)^{-\frac{n}{2}}\int_M \tau R_{ij}h_{ij}
e^{-f}\mrm{d}V.
\end{align*}
Therefore,
$$\delta\tau = \tau \frac{\int_M R_{ij}h_{ij} e^{-f}\mrm{d}V}{\int_M R e^{-f}}.$$ 
\end{proof}

Now,  by Lemmas \ref{lem2.3} and  \ref{lem2.4}, the second variation becomes

\begin{align*}
 \delta^2 \nu_g (h, h)= & (4\pi\tau)^{-\frac{n}{2}}\int_M -\tau <h, \delta Rc+\delta\nabla^2 f-\frac{1}{2\tau} h> e^{-f} \\
& + (4\pi\tau)^{-\frac{n}{2}}(-\frac{\delta\tau}{2\tau}) \int_M \tr_g h  e^{-f}\\
=&  (4\pi\tau)^{-\frac{n}{2}}\int_M \tau <h, \frac{1}{2} \Delta_f h+Rm (h, \cdot) +\div_f^{\dagger}\div_f h>e^{-f} \\
& + (4\pi\tau)^{-\frac{n}{2}}\int_M \tau <h, \nabla^2 (-\delta f +\frac{1}{2} \tr_g h)> e^{-f} \\
& + (4\pi\tau)^{-\frac{n}{2}}(-\frac{\delta\tau}{2\tau}) \int_M \tr_g h  e^{-f}\\
= &  \tau (4\pi\tau)^{-\frac{n}{2}}\int_M <h, \frac{1}{2} \Delta_f h+Rm (h, \cdot) +\div_f^{\dagger}\div_f h+\frac{1}{2}\nabla^2 \hat v_h>e^{-f} \\
& + \tau (4\pi\tau)^{-\frac{n}{2}}\frac{\delta\tau}{\tau} \int_M <h,\nabla^2 f-\frac{1}{2\tau} g>  e^{-f}\\
= & \tau (4\pi\tau)^{-\frac{n}{2}}\int_M <h, \frac{1}{2} \Delta_f h+Rm (h, \cdot) +\div_f^{\dagger}\div_f h+\frac{1}{2}\nabla^2 \hat v_h>e^{-f} \\
& - \tau (4\pi\tau)^{-\frac{n}{2}}  \frac{\int_M <Rc, h>e^{-f}dV}
{\int_M R e^{-f}dV} \int_M <h,Rc>  e^{-f}dV.
\end{align*}
Here, $$\hat v_h=-2\delta f +\tr_g h-\frac{2\delta\tau}{\tau}
(f-\nu),$$ and it is straightforward to check that
$$\Delta_f \hat v_h+\frac{\hat v_h}{2\tau}=\div_f\div_f
h,\qquad\int_M \hat v_h e^{-f}dV=0. \eqno(2.15)$$

To see the uniqueness of the solution to (2.15), it suffices to show that $\lambda_1 (\Delta_f)>\frac{1} {2\tau}$, where $\lambda_1=\lambda_1 (\Delta_f)$ denotes the first eigenvalue
of $\Delta_f$.  Let $u$ be a (non-constant)
first eigenfunction so that $$ \Delta_f u=-\lambda_1 u. $$ Then by direct computations (see also \cite{Futaki}), we get
\begin{align*}
\frac 1 2 \Delta_f |\nabla u|^2= &|\nabla^2 u|^2 +\nabla (\Delta_f u)\cdot \nabla u +(Rc+\nabla^2 f) (\nabla u,\nabla u)\\
\ge  & \frac 1 n |\Delta u|^2  +(\frac {1} {2\tau}-\lambda_1)|\nabla u|^2.
\end{align*}
Thus,
$$0=\int_M \frac 1 2 \Delta_f |\nabla u|^2 e^{-f} dV\ge  \frac 1 n \int_M |\Delta u|^2e^{-f}dV   +(\frac {1} {2\tau}-\lambda_1)\int_M |\nabla u|^2e^{-f} dV.$$
Since $u$ is non-constant, we obtain   $$\lambda_1>\frac {1} {2\tau}.$$

This completes the proof of Theorem \ref{thm1.1}.

\qed

\section{Further remarks and the proof of Theorem \ref{thm1.3}} \label{sec:3}

Recall that a gradient shrinking Ricci soliton $(M^n, g_{ij}, f)$
is called {\it linearly stable} if the stability operator $\hat
N\le 0$ on symmetric 2-tensors. Note that $\hat N$ is degenerate
negative elliptic. In this section we shall exhibit the action of
the stability operator on a couple of special symmetric 2-tensors
$h$: (i) $h_{ij}=g_{ij}$ and (ii) $h_{ij}=R_{ij}$, and prove
Theorem 1.3.

Without loss of generality, we assume $\tau=1$ so that our
shrinking soliton $(M^n, g_{ij}, f)$ satisfies the equation
$$
R_{ij}+\nabla_i\nabla_j f=\frac{1}{2}g_{ij}. \eqno(3.1)
$$
We also normalize $f$ so that
$$(4\pi)^{-\frac{n}{2}}\int_M e^{-f}=1.$$

First of all, notice that we have
$$\Delta_f g=\Delta g-\nabla f \nabla g=0$$
$$\div_f g=-g_{ij}\nabla_j f=-\nabla f, \qquad \div_f^{\dagger} \div_f g=\nabla^2 f,$$
$$\hat v_h=-2(f-\bar f), \qquad \bar f=(4\pi\tau)^{-\frac{n}{2} }\int_M f e^{-f}.$$
Hence, we get
$$
\hat N (g)=Rc +\nabla^2 f +\frac{1}{2}\nabla^2 (-2f+2\bar f) -Rc=0  \eqno(3.2)
$$
as we expected.  \\

On the other hand, we have

\begin{lemma} \label{lem3.1} For any complete shrinking Ricci soliton satisfying (3.1), we have $$ Rc \in \ker\div_f.$$
\end{lemma}

\begin{proof} By definition and the second contracted Bianchi identity, $$(\div_f Rc)_i=\nabla_j R_{ij}-R_{ij}\nabla_jf =\frac 12 \nabla_i R-R_{ij}\nabla_jf.$$
On the other hand, it is a basic fact that our shrinker satisfies
$$
\nabla_i R=2R_{ij}\nabla_jf.  \eqno(3.3)
$$
Therefore,  $\div_f (Rc)=0$.

\end{proof}

Recall the operator $\mathcal{L}_f$ on symmetric 2-tensors defined
in (1.10):

$$
\mathcal{L}_f h:=\frac{1}{2}\Delta_f h +Rm(h,\cdot ).
$$
It is easy to see that $\mathcal{L}_f $ is a self-adjoint operator
with respect to the weighted $L^2$-inner product $(\cdot,
\cdot)_f$ defined in (1.6).

\begin{lemma}\label{lem3.2}  For any complete shrinking soliton satisfying (3.1), its Ricci tensor is an eigen-tensor of the operator $\mathcal{L}_f$ with eigenvalue
$-1/2$:
$$ \mathcal{L}_f (Rc)=\frac 12 Rc.$$
\end{lemma}

\begin{proof} The following computations are familiar to experts, but we carry out here for completeness.

From the soliton equation (3.1), we have $$R_{ij}=\frac 1 2
g_{ij}-\nabla_i\nabla_j f.$$ By commuting covariant derivatives,
we have
\begin{align*}
\Delta R_{ij} &= -\nabla_k\nabla_k\nabla_i\nabla_j f\\
&=-\nabla_k(\nabla_i\nabla_k\nabla_j f +R_{kijl}\nabla_lf)\\
&= -\nabla_k\nabla_i\nabla_j\nabla_k f - \nabla_kR_{kijl}\nabla_lf
-R_{kijl}\nabla_k\nabla_lf.
\end{align*}

On the other hand, by commuting covariant derivatives again and
using the contracted second Bianchi identity as well as (3.1), we
obtain
\begin{align*}
 \nabla_k\nabla_i\nabla_j\nabla_k f & =\nabla_i\nabla_k\nabla_j\nabla_k f +R_{kijl}\nabla_l\nabla_kf +R_{kikl}\nabla_j\nabla_lf \\
 &=-\frac 1 2 \nabla_i\nabla_j R-R_{kijl}R_{kl} -R_{il}R_{lj}\\
 &=-\nabla_jR_{il} \nabla_lf-\frac 12 R_{ij} -R_{kijl}R_{kl}.
 \end{align*}
 Here we have used (3.3) in deriving the last equality.

Moreover,  by the second Bianchi identity, we have
$$\nabla_kR_{kijl}\nabla_lf=(\nabla_jR_{il}-\nabla_lR_{ij})\nabla_lf.$$

Combining the above calculations and using the Ricci soliton
equation (3.1), we arrive at
$$
\Delta R_{ij} =\nabla_lR_{ij}\nabla_lf +2 R_{kijl}R_{kl} +R_{ij},  \eqno(3.4)
$$
i.e., $2\mathcal{L}_f (R_{ij})=R_{ij}$. 

\end{proof}

Now, for any $h\in \ker\div_f$, the stability operator $\hat N$ is
given by
$$
\hat N h=\mathcal{L}_f h -Rc\: \frac{\int_M <Rc, h>e^{-f}}{\int_M Re^{-f}}.  \eqno(3.5)
$$
Moreover, from (3.4) we obtain
$$
\Delta_f R=R-2|Rc|^2,  \eqno(3.6)
$$
from which it follows that
$$
2\int_M |Rc|^2 e^{-f}= \int_M Re^{-f}.  \eqno(3.7)
$$
Therefore, by Lemma \ref{lem3.2}, (3.5) and (3.7), we have
\begin{lemma} \label{lem3.3}
$$\hat N (Rc)=0.$$
\end{lemma}

Now we are ready to prove

\begin{proposition} \label{prop3.1}
Suppose $(M^n, g_{ij}, f)$ is a linearly stable compact shrinking soliton satisfying (3.1),
then $-1/2$ is the only negative eigenvalue of the operator
$\mathcal{L}_f $ on $\ker\div_f$, and the multiplicity of  $-1/2$
is one.  In particular, $-1/2$ is the least eigenvalue of
$\mathcal{L}_f $ on $\ker\div_f$.
\end{proposition}

\begin{proof} By Lemma \ref{lem3.1} and Lemma \ref{lem3.2}, we know that $Rc\in \ker\div_f$, and is an eigen-tensor
of $\mathcal{L}_f $ with eigenvalue $-1/2$. Suppose there exists a
(non-zero) symmetric 2-tensor $h\in \ker\div_f$ such that
$$\mathcal{L}_f h =\alpha h,$$ with $\alpha >0$, and $$(Rc, h)_f=:\int_M <Rc, h> e^{-f}=0.$$
Then, by Theorem \ref{thm1.1} and (3.5), we have
\begin{align*}
\delta^2\nu_g (h,h) =& \frac{1}{(4\pi)^{n/2}}\int_M  <\hat{N}h, h> e^{-f} \\
=& \frac{1}{(4\pi)^{n/2}}\int_M  <\mathcal{L}_f h, h> e^{-f}\\
=&  \frac{\alpha}{(4\pi)^{n/2}}\int_M  |h|^2 e^{-f}>0,
\end{align*}
a contradiction to the linear stability of $(M^n, g_{ij}, f)$.
Thus $-1/2$ is the only negative eigenvalue of $\mathcal{L}_f $ on
$\ker\div_f$, with multiplicity one.

\end{proof}

\begin{remark} \label{rmk3.1} In \cite{HM}, the authors have given a very nice interpretation of their proof
in terms of the multiplicity of the eigenvalue $-1/2$: for any
compact shrinking K\"ahler-Ricci soliton satisfying (3.1), the
eigen-space of eigenvalue $-1/2$ has multiplicity at least
$h^{1,1}$. Hence a compact shrinking K\"ahler-Ricci soliton with
$h^{1,1}>1$ is unstable.

\end{remark}

\end{document}